\documentclass[11pt]{article}
\usepackage{amsfonts}
\usepackage{amsmath}
\usepackage{amssymb}
\usepackage{amsthm}
\usepackage{amsbsy}
\usepackage{paralist}
\usepackage{graphicx}
\usepackage{tikz}
\usepackage{filecontents,pgfplots}
\usepackage{comment}
\usepackage[colorlinks=true]{hyperref}
\hypersetup{urlcolor=blue, citecolor=red}
\usepackage[width=0.9\textwidth,font=small]{caption}

\newcommand{\rest}[2]{{#1}\lower0.5ex\hbox{$|$}{}_{#2}}

\renewcommand{\aa}{{\boldsymbol a}}
\newcommand{\eps}{{\varepsilon}}

\newcommand{\PD}[2]{\frac{\partial{#1}}{\partial{#2}}}

\newcommand{\R}{{\mathbb R}}

\newcommand{\calA}{\mathcal{A}}

\newcommand{\calI}{\mathcal{I}}
\newcommand{\calJ}{\mathcal{J}}

\newcommand{\calO}{{\mathcal O}}

\newcommand{\calU}{{\mathcal U}}

\newcommand{\hatO}{{\widehat\Omega}}

\renewcommand{\AA}{{\boldsymbol A}}

\newcommand{\KK}{{\boldsymbol K}}

\newcommand{\bb}{{\boldsymbol b}}

\newcommand{\dd}{{\boldsymbol d}}
\newcommand{\ff}{{\boldsymbol f}}
\newcommand{\qq}{{\boldsymbol q}}

\newcommand{\uu}{{\boldsymbol u}}

\newcommand{\yy}{{\boldsymbol y}}

\newcommand{\xx}{{\boldsymbol x}}

\newcommand{\tran}[1]{{#1}^{\text{T}}}

\newcommand{\minimize}{\operatornamewithlimits{minimize}}

\theoremstyle{definition}

\newtheorem{example}{Example}

\begin{document}

\title{On optimal control of free boundary\\ problems of obstacle type} 
\author{Raino A. E. M\"akinen\thanks{Faculty of Information
    Technology, University of Jyvaskyla, Finland. E-mail: {\tt
      raino.a.e.makinen\,@\,jyu.fi}}} 
\date{\today}
\maketitle

\abstract{A numerical study of an optimal
  control formulation for a shape optimization 
  problem governed by an elliptic variational inequality is performed. The shape
  optimization problem is reformulated as a boundary control problem
  in a fixed domain. The
  discretized optimal control problem
  is a non-smooth and non-convex mathematical programing problem. 
  The performance of the standard BFGS quasi-Newton method and
  the BFGS method 
  with the inexact line search are tested.}

\medskip



\section{Introduction}

In the present paper a numerical method for shape optimization for a
class of free boundary problems is proposed. The shape optimization
problem is reformulated as a boundary control problem following the
ideas discussed in \cite{NT95}.
The optimal shape is obtained as a level set. Thus, in numerical
realization using finite elements, there is no need for remeshing
during optimization iterations. 

As a concrete free boundary (state) problem we consider the contact problem
for the Poisson problem. In this case the free boundary is defined by
the a-priori unknown contact zone. As  a shape optimization  problem
we consider the problem of finding a shape having minimum area such
that the contact zone includes a given subdomain \cite{BSZ84}.

It is well-known that in control problems for contact problems the
control-to-state mapping is not smooth, in general. Despite this, many
nonlinear programming codes for smooth optimization can be used to
obtain reasonable approximate solutions to those control problems. 

The paper is organized as follows. In Section 2 we present a shape
optimization problem and an alternate boundary control formulation of
it. Section 3 is devoted to disceretization and algebraic sensitivity analysis.
Finally, 
a numerical study on the performance of two nonlinear programming codes
applied to the discretized control problem is done in Section~4.

\section{Setting of the problem}

We consider an abstract shape optimization problem given formally as follows: 
\begin{equation}
\left\{\begin{aligned}
   &\minimize_\alpha F(\alpha,\Omega_\alpha,y_\alpha) \\
\text{subject to}&\\
   &\text{Find }(\Omega_\alpha,y_\alpha): \quad
\calA(\alpha,\Omega_\alpha,y_\alpha)=0. 
\end{aligned}\right.
   \label{eq:abs-opt}
\end{equation}
Here $\alpha$ is the optimization (control) parameter, $\calA$ represents a {\em
free boundary} problem (in abstract operator form) whose solution
consists both of the domain $\Omega_\alpha$  where the PDE is posed as well
as the solution $y_\alpha$ to that PDE. 

Problem \eqref{eq:abs-opt} is more complicated than a classical shape
optimization problem. In the latter the state problem is solved on a fixed
domain, in the sense that $\Omega_\alpha$ is considered to be given
instead of being one of the unknowns in the state problem.


Next we introduce the model problem to be studied in the rest of this paper.
Let $\Omega\subset\R^2$ a domain, $f\in L^2(\Omega)$, $\psi\in C^1(\Omega)$.
We consider the following classical free boundary problem:
\medskip
Find $(y_\Omega,Z_\Omega),\ Z_\Omega{\subset}\Omega$ such that
\begin{equation}
  \left\{\begin{aligned}
    -\Delta y_\Omega=f,\ y_\Omega>\psi &\quad\text{in }\Omega\setminus Z_\Omega\\
      y_\Omega =0& \quad\text{on }\partial\Omega\\
      y_\Omega=\psi &\quad\text{in }Z_\Omega\\
      \PD{y_\Omega}{n}+\PD\psi n=0 &\quad\text{on }\partial Z_\Omega.
  \end{aligned}\right.
    \label{eq:FB0}
\end{equation}
Problem \eqref{eq:FB0} models e.g.\ the contact between elastic
membrane (represented by $\Omega$) and a rigid obstacle defined by
$\psi$. The solution $y_\Omega$ to  \eqref{eq:FB0} then gives the
vertical displacement of the membrane under the vertical load
represented by $f$. The set $Z_\Omega$ is the contact zone between the
membrane and the obstacle.

Problem \eqref{eq:FB0} can be formulated as the following variational
inequality (\cite{Glow84}):
\begin{equation}
  y_\Omega\in K_0(\Omega): \quad \int_\Omega  \nabla y_\Omega 
\cdot  \nabla(w{-}y_\Omega)\,dx\ge\int_\Omega f(w-y_\Omega)\,dx \quad
w\in K_0(\Omega),
  \label{eq:VI}
\end{equation}
where
$$
K_g(D):=\left\{ v\in H^1(D) \mid   v{\ge}\psi\text{ in } D\text{ and }
v{=}g \text{ on }\partial D\right\}.
$$ 
This formulation does not include explicitly the contact zone
$Z_\Omega$ as an unknown.

Let $\hatO\subset\R^2$ and $\Omega_0\subset\R^2$  be given bounded domains
such that $\Omega_0\subsetneq\hatO$. We introduce  a system of
bounded domains
$$
\calO=\{ \Omega\subset\hatO \mid \Omega \text{ has Lipschitz
  boundary}\}.
$$

Let $y_\Omega$ be the unique solution to \eqref{eq:VI} and 
let $C(y_\Omega)=\{ x\in\Omega \mid y_\Omega(x)=\psi(x)\text{ a.a. }
x\}$ denote the 
corresponding contact region. 
We consider the following shape optimization problem introduced in
\cite{BSZ84}: 
\begin{equation}
\left\{\begin{aligned}
  &\minimize_{\Omega\in\calO}\ \calJ(\Omega)=\int_\Omega\,dx \\
   &\text{subject to \eqref{eq:VI} and } \Omega_0\subset C(y_\Omega).
\end{aligned}\right.
  \label{eq:PB-Ome}
\end{equation}
We are thus looking for a domain $\Omega^\star$ (representing the
membrane) having minimum area
such that the contact region ''covers'' the given domain $\Omega_0$.

The state constraint $\Omega_0\subset C(y_\Omega)$ can be relaxed by
introducing a penalty term resulting in the shape optimization problem
\begin{equation}
\left\{\begin{aligned}
  &\minimize_{\Omega\in\calO}\quad \calJ_\eps(\Omega)=\int_\Omega\,dx +
  \frac1{\eps}\int_{\Omega_0}(y_\Omega-\psi)^2\,dx\\ 
   &\text{subject to \eqref{eq:VI}}, \ \eps>0.
\end{aligned}\right.
  \label{eq:PB-Ome-eps}
\end{equation}




One of the main difficulties in the theoretical and numerical
treatment of the problem 
\eqref{eq:PB-Ome-eps} 
consists in the variable character of the domain $\Omega$ on which the
state problem \eqref{eq:VI} is given.  
Theoretical and numerical aspects of shape optimization using
''moving'' domains are discussed e.g.\ in refs 
\cite{Pir84}, \cite{Sokolowski92}, \cite{HN96}, \cite{has03}. 

Alternatively, there has been significant amount of interest in fixed
domain formulations of shape optimization problems
(see e.g.\ \cite{Fulmanski07}). 
The most common
approaches to get rid of moving domains are the following:
\begin{itemize}
\item By scaling 
the domain such that it becomes fixed, the optimization
parameter appears as 
a coefficient in the state problem. 
\item The domain $\Omega$ is represented by a level set
  $\Omega=\{x\in\hatO\mid \Phi(x)<0\}$, where $\Phi:\hatO\to\R$ is an
  unknown level set function to be determined.
\item The state problem is modified by adding control variable to
  it. An optimal shape is defined implicitly by a level set (but
  without separate level set function).
\end{itemize}

\subsubsection*{Boundary control approach}

In this paper we utilize the boundary control approach discussed in \cite{NT95}.
In what follows we assume that $f\le0$ and $\psi<0$.
We define the set of admissible controls as follows
$$
\calU=\{ u:\partial\hatO\to\R \mid u\in C(\partial\hatO),
\ u_{\min}\le u\le u_{\max}\}, 
$$
where $u_{\min},u_{\max}>0$ are given constants.
Next, to each $u\in\calU$ we associate $y(u)$,  the solution to the
variational inequality 
\begin{equation}
\left\{\begin{aligned}
 &\text{Find } y(u)\in K_u(\hatO) \text{ such that} \\
&\int_\hatO \nabla y(u)\cdot\nabla(w{-}y(u))\,dx 
\ge \int_\hatO f(w{-}y(u))\,dx \quad \forall
w\in K_u(\hatO).
\end{aligned}\right.
\label{eq:VI2}
\end{equation}

Let $H:\R\to\R$ denote the Heaviside step function. We consider the
following boundary control problem posed in a fixed domain:
\begin{equation}
 \left\{\begin{aligned}
 &\minimize_{u\in\calU}\ 
  J_\eps(u) = \int_{\hatO} \bigl(1-H(y(u))\bigr) \,dx + \frac1\eps
  \int_{\Omega_0}(y(u)-\psi)^2\,dx \\
 &\text{subject to \eqref{eq:VI2}}.
 \end{aligned}\right.
  \label{eq:PB-u}
\end{equation}
As $u>0$ on $\partial\hatO$ and $f\le0,\psi<0$ in $\hatO$ it follows
that for suitable choice of $u$ the set of points where $y(u)$ is
strictly negative is non-empty. The first term in the cost functional
$J_\eps$ gives the area of this set while the latter term adds  a
penalty if the contact zone $C(y(u))$ does not cover $\Omega_0$.

\def\ues{{u_\eps^\star}}

Assume that there exists an optimal pair
$(u^\star_\eps,y^\star_\eps):=(\ues,y(\ues))$ for
\eqref{eq:PB-u}. We can now 
consider $\Omega_\eps^\star := \{x\in\hatO\mid y^\star_\eps(x)<0\}$ as an
approximate solution to the original shape optimization problem
\eqref{eq:PB-Ome}. 

The control approach clearly makes sense if $\Omega_0$ is simple
enough star-like domain. If, however, $\Omega_0$ is e.g.\ multiply
connected the set $\Omega_\eps^\star$ might not approximately solve
\eqref{eq:PB-Ome-eps}.

\section{Approximation and numerical realization}

To realize \eqref{eq:PB-u} numerically, we must discretize both the
control and state variables. In what follows, we assume that $\hatO$ is
the disk $B(0,R)$. We discretize the state problem
\eqref{eq:VI2} by using piecewise linear triangular elements. 
Instead of the exact Heaviside function, we use the smoothed one
$$
    H_\beta(x) := \tfrac12\tanh(\tfrac x\beta)+\tfrac12, \quad \beta>0.
$$

The use of piecewise linear triangular elements implies that the
most obvious way to discretetize the control $u$ would be to use piecewise
linear and continuous discretization in the same finite element mesh
where the state variable is discretized. However, this approach has
two well-known drawbacks. Firstly, the number of optimization variables is very
large whenever dense meshes are used. Secondly, the piecewise linear
approximation (without a suitable regularization term in the cost
function) is prone to spurious oscillations. 

Instead, we look for a differentiable and periodic function
$u_\aa:[0,2\pi]\to\R$ 
that is fully defined by a vector of parameters  $\aa=(a_1,a_2,...,a_n)$.
Examples of such functions are Bezier functions, cubic Hermite, and cubic spline
interpolating polynomials, e.g.   

Here we shall use shape preserving periodic piecewise cubic Hermite
interpolation 
polynomial (see \cite{Fritsch80}, \cite{Bastien14}) to represent the
control function. The parametrized control function solves the
following interpolation problem
\begin{align*}
 &u_\aa\in C^1([0,2\pi]), \quad u'_\aa(0-)=u_\aa'(2\pi+)\\
 &u_\aa(0){=}a_1,\ u_\aa(\Delta){=}a_2,...,u_\aa(2\pi{-}\Delta){=}a_n,\ 
 u_\aa(2\pi){=}a_1.
\end{align*}
The advantage of this kind of parametrization is
that the number of discrete optimization variables is small but at the
same time the control is a smooth function without excessively prone
to wild oscillations. Moreover, the interpolant $u_\aa$ does not
overshoot the data, so $u_\aa\in\calU$ if $u_{\min}\le a_i \le
u_{\max}\ \forall i$.

\subsubsection*{Sensitivity analysis for the discrete state problem}

Let $\calI=\{1,...,N\}$ denote the set of node numbers in the finite
element model and let $\calI_d\subset\calI$ be the node numbers of the
boundary nodes. Let $\boldsymbol\psi=\{\psi_1,...,\psi_N\}$ be the
vector of
nodal values of the obstacle and let $\uu=\{u_{i_1},...,u_{i_m}\}$ be
the vector of
nodal values of the boundary control.
Then the finite element approximation of the variational inequality
\eqref{eq:VI2}  can be expressed as an equivalent quadratic programming problem
\begin{equation}
 \left\{\begin{aligned}
   &\minimize_{\qq}\ \Pi(\qq)=\tfrac12 \tran\qq\KK\qq-\tran\ff\qq \\
   &\text{subject to}\quad q_i{=}u_\aa(x_i), \ i\in \calI_d, \quad
   q_i{\ge}\psi_i, \ i\in\calI\setminus\calI_d,
 \end{aligned}\right.
 \label{eq:QP}
\end{equation}
where $\KK$ and $\ff$ are the stiffness matrix and the force vector,
respectively. 

\def\taa{\widetilde{\aa}}

The discretization of the boundary control problem \eqref{eq:PB-u}
leads to the nonlinear programming problem
\begin{equation}
   \left\{\begin{aligned}
   &\minimize_{\aa} J_\eps(\aa) \\
   &\text{subject to } \eqref{eq:QP} \text{ and } u_{\min}\le a_k\le
   u_{\max},\ k=1,...,n.
   \end{aligned}\right.
\end{equation}
To be able to use descent type optimization methods we need to
evaluate the gradient of $J_\eps$ with respect to the discrete control
variable $\aa$.

Let $\calI_c(\aa)$ denote the set of contact nodes, i.e.\ the solution
$\qq=\qq(\aa)$ satisfies
$q_i=\psi_i,\ i\in\calI_c(\aa)$. Let $\taa$ be the current
approximation of the optimal control.
If we assume that $\calI_c(\taa)$ is known a
priori and it is invariant under small perturbations of $\taa$, then
$J$ is differentiable and the gradient $\nabla_\aa J_\eps(\taa)$ can
be obtained using the standard adjoint equation technique. However, if
the contact set changes due to arbitrary small perturbation, the
mapping $\aa\mapsto J_\eps(\aa)$ is not differentiable at $\taa$ and
the mechanical application of the adjoint technique gives at most an
element from the subdifferential.

\subsubsection*{Additional source of nonshmoothness}

Consider piecewise cubic Hermite interpolation of the
data $(x_1,y_1)$, $(x_2,y_2),....,(x_n,y_n)$. On each subinterval
$[x_k,x_{k+1}]$ the interpolant can be expressed in terms of a local
variable $s=x-x_k$ as follows
\begin{equation*}
 p(x)=\frac{3hs^2-2s^3}{h^3}y_{k+1} + \frac{h^3-3hs^2+2s^3}{h^3}y_k +
 \frac{s^2(s-h)}{h^2}d_{k+1} + \frac{s(s-h)^2}{h^2}d_k.
\end{equation*}
The interpolant satisfies the following interpolation conditions
$$
p(x_k)=y_k,\ p'(x_k)=d_k,\ p(x_{k+1})=y_{k+1}, \ p'(x_{k+1})=d_{k+1}.
$$
In the classical piecewise cubic Hermite interpolation, the slope parameters
$d_k,\ k=1,...,n$ are a priori given constants.

In shape-preserving piecewise cubic interpolation, the slopes are
not supplied by the user. Instead, they are computed algorithmically 
as a part of the interpolating process in the following way:
Let $\delta_k=(y_{k+1}-y_k)/(x_{k+1}-x_k)$ and $h_k=x_{k+1}-x_k$. The
slopes $d_k$  are  
determined as follows.
If $\text{sign}(\delta_k)\cdot\text{sign}(\delta_{k-1})<0$ we set
  $d_k=0$. 
If $\text{sign}(\delta_k)\cdot\text{sign}(\delta_{k-1})>0$
    and the two intervals have the same length, then $d_k$ is the
    harmonic mean of $\delta_{k-1}$ and $\delta_k$:
$$
   \frac1{\delta_k} = \tfrac12\left(
   \frac1{\delta_{k-1}}+\frac1{\delta_k}\right). 
$$
If the intervals have different lengths, then $d_k$ is the following
weighted harmonic mean
$$
  \frac{w_1+w_2}{d_k}=\frac{w_1}{\delta_{k-1}}+ \frac{w_2}{\delta_k},
$$
where $w_1=2h_k+h_{k-1}$, $w_2=h_k+2h_{k-1}$.

The advantage of shape-preserving piecewise cubic interpolation over
cubic spline interpolation is
that the interpolant does not overshoot the data. 

Consider next the case where the vector of values $\yy=(y_1,...,y_n)$ is
an optimization 
parameter. In case of cubic spline interpolation, the vector of slopes
is obtained by solving a tridiagonal system  $\AA\dd=\bb$, where the
vector $\bb$ is a linear function of $\yy$. Thus the value of the
interpolant at given fixed point $x$ is a smooth function of $\yy$. 

\begin{figure}[h]
\begin{center}
\includegraphics[scale=0.65]{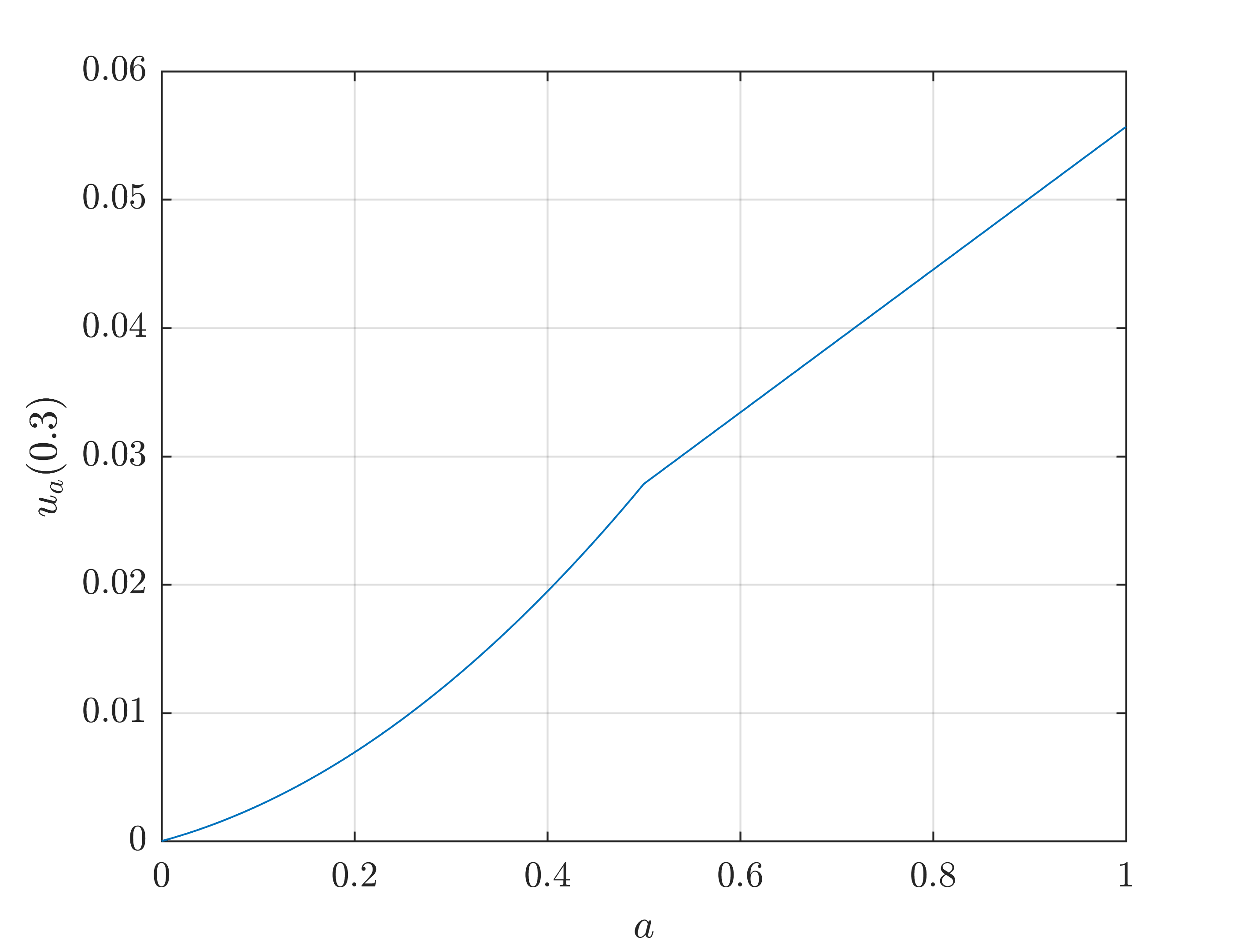}
\end{center}
\caption{Periodic shape-preserving piecewise cubic Hermite interpolation of
  the data
  $\xx{=}(0,\tfrac{2\pi}3,\tfrac{4\pi}{3},2\pi),\ 
  \yy{=}(0, a, \tfrac12, 0)$. The value of the interpolant $u_a(x)$ at
  $x{=}0.3$ is plotted against the parameter $a$.
}\label{fig:pchip}
\end{figure}

With shape-preserving piecewise cubic interpolation the situation is
different. The mapping $\yy\mapsto\dd$ is clearly nonlinear. In
addition, it is {\em nonsmooth}. A simple example shown in
Figure~\ref{fig:pchip} demonstrates this.

\section{Numerical examples}

In this section we present two numerical examples. All finite element
computations are implemented in  Matlab \cite{matlab16}.
The quadratic programming problem \eqref{eq:QP} is solved using Matlab's {\tt
  quadprog} function. 
The box constraints $u_{\min}\le a_i\le u_{\max}$ are implemented by
adding a smooth penalty function 
$$
  \tfrac\gamma2\sum_{i=1}^n\left[ (a_i-u_{\max})_+^2 +
    (u_{\min}-a_i)_+^2\right]
$$
with $\gamma=10^{-3}$ to the cost function.

\goodbreak
The optimal control problem is solved by using
two different optimizers: 
\begin{itemize}
\item Function {\tt fminunc} from Matlab's Optimization Toolbox with
  {\tt 'quasi-newton'} option that implements the standard BFGS method.
\item Optimizer HANSO~2.2 (Hybrid Algorithm for Non-Smooth Optimization
  \cite{hanso1}, \cite{hanso2}) that implements BFGS method 
with inexact line search. 
\end{itemize}

For the options of {\tt fminunc} that control the optimality tolerance
we use their default values. 
Also we use the default options of HANSO. Especially, we do not use
the option available in HANSO to continue the 
optimization using the gradient sampling method as this option is very
expensive in terms of the number of function evaluations.

In all examples we use the following parameter values. The
radius defining  the computational domain $\hat\Omega$ is $R=1.75$.
The right hand side function $f\equiv-10$ and the obstacle 
$\psi(x)= -0.3\left(x_1^2 + (x_2{-}0.25)^2\right)-0.05$. The target
domain $\Omega_0$ to be covered is the isosceles triangle with
vertices $({-}1,0)$, 
$(\tfrac12,\tfrac34)$, $(\tfrac12,{-}\tfrac32)$
The parameters defining $\calU$ are 
$u_{\min}{=}0.01$ and $u_{\max}{=}10$. The Heaviside smoothing parameter
is $\beta{=}10^{-3}$ and penalty parameter is $\eps=10^{-3}$.

\begin{example}\label{ex:eka}
Here we use an unstructured finite element mesh with nominal mesh size
$h=0.05$. The mesh is constructed in such way that $\partial\Omega_0$
coincides element edges.


We ran both optimizers using $n=30$ and
$a_i^{(0)}{=}2,\ i{=}1,...,n$ until they were not able to improve the control.
The evolution of the best objective value versus computational work is
depicted for both methods in Figure~\ref{fig:h005cost}. 
The implicitly defined optimized domains $\Omega_\eps^\star$ (with the
corresponding contours of $y_\eps^\star$) and the optimal boundary
controls are plotted in Figures~\ref{fig:h005} and \ref{fig:h005H}.

\begin{figure}[h]
\begin{center}
 \includegraphics[scale=0.6]{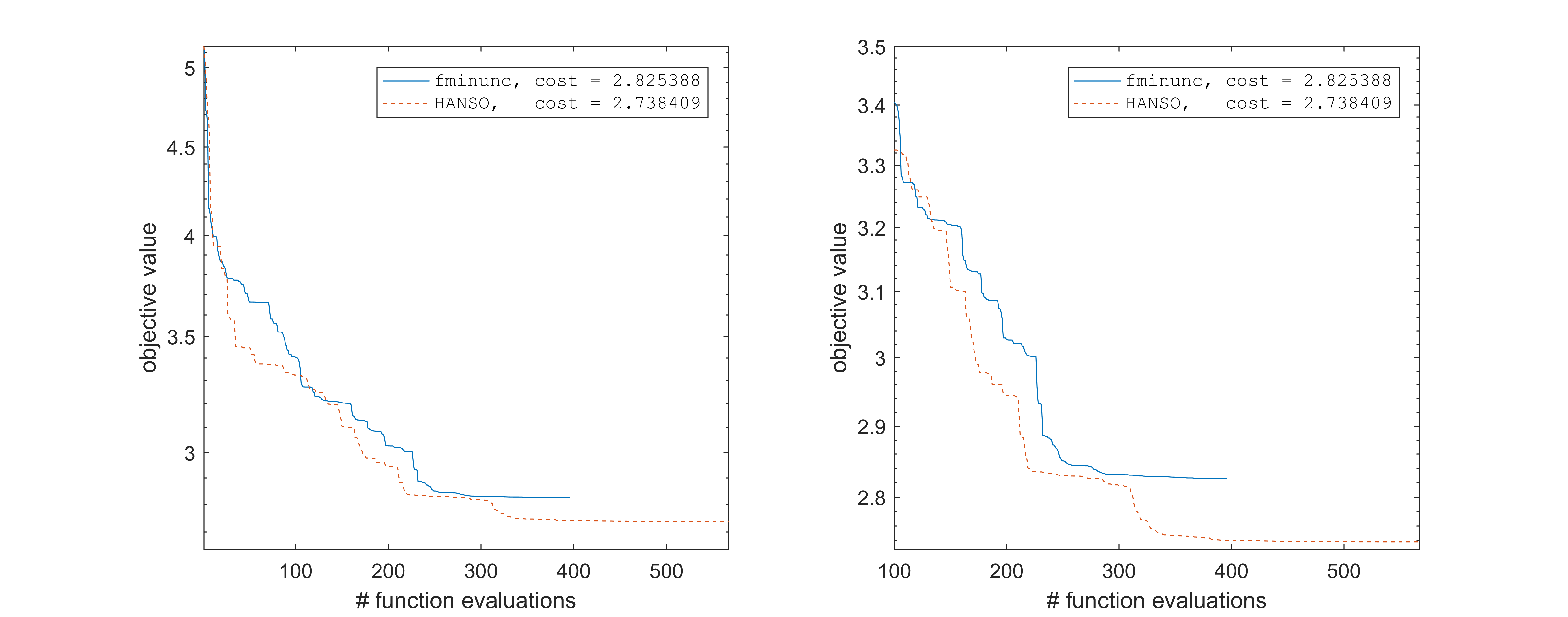}
\end{center}
\caption{Best objective value versus the number of function
  evaluations}\label{fig:h005cost} 
\end{figure}

\begin{figure}[h]
\begin{center}
 \includegraphics[scale=0.6]{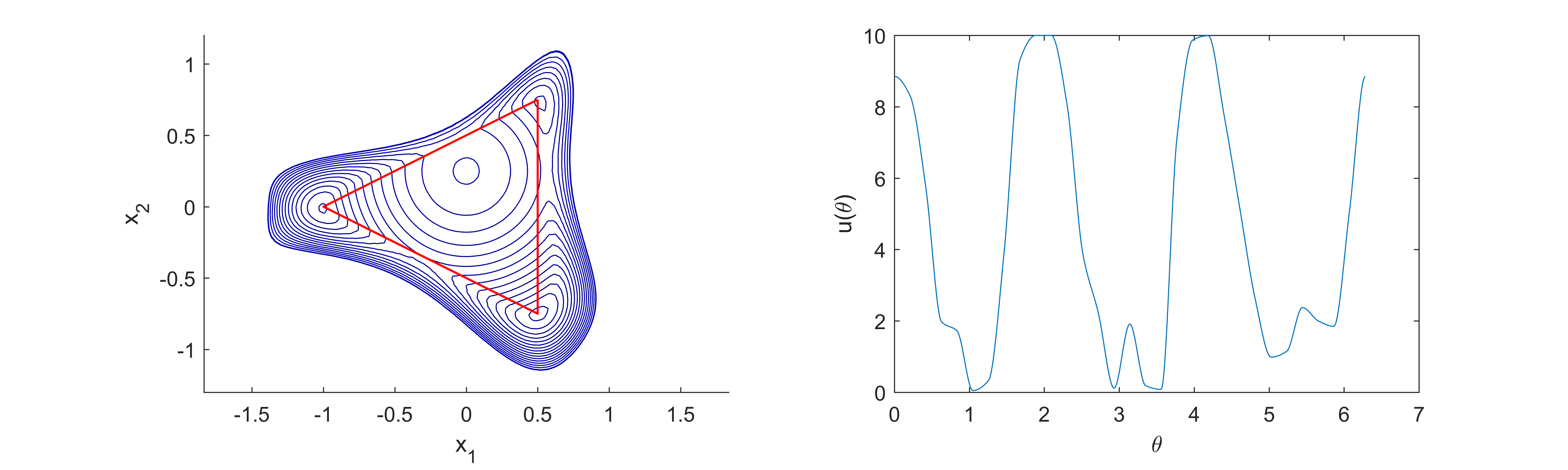}
\end{center}
\caption{Optimized domain (left) and control (right) computed using
  fminunc. The boundary of the target contact zone $\Omega_0$ is
  depicted as the red triangle.}\label{fig:h005} 
\end{figure}

\begin{figure}[h]
\begin{center}
 \includegraphics[scale=0.6]{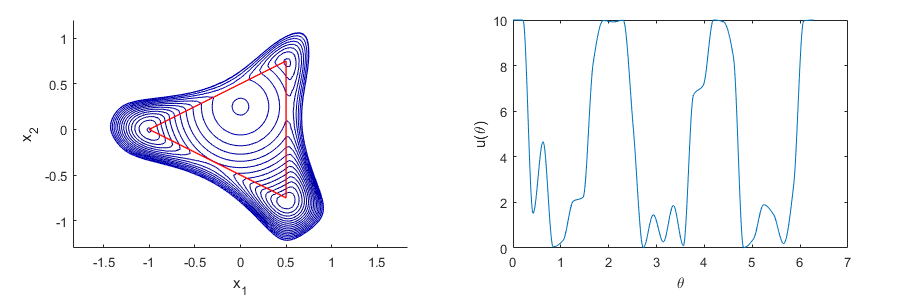}
\end{center}
\caption{Optimized domain (left) and control (right) computed using
  HANSO~2.2}\label{fig:h005H} 
\end{figure}

One can observe that in this example HANSO~2.2 finds 
better control than {\tt fminunc}. However, HANSO~2.2 spends
considerable amount of computational work with only very marginal
improvement in the cost function.

\end{example}

\begin{example}
In this example the mesh is generated in the same way as in
Example~\ref{ex:eka} except the nominal mesh size is $h=0.025$.
The evolution of the best objective value versus computational work is
depicted for both methods in Figure~\ref{fig:h0025cost}. 
The implicitly defined optimized domains $\Omega_\eps^\star$ (with the
corresponding contours of $y_\eps^\star$) and the optimal boundary
controls are plotted in Figures~\ref{fig:h0025} and \ref{fig:h0025H}.

\begin{figure}[h]
\begin{center}
 \includegraphics[scale=0.6]{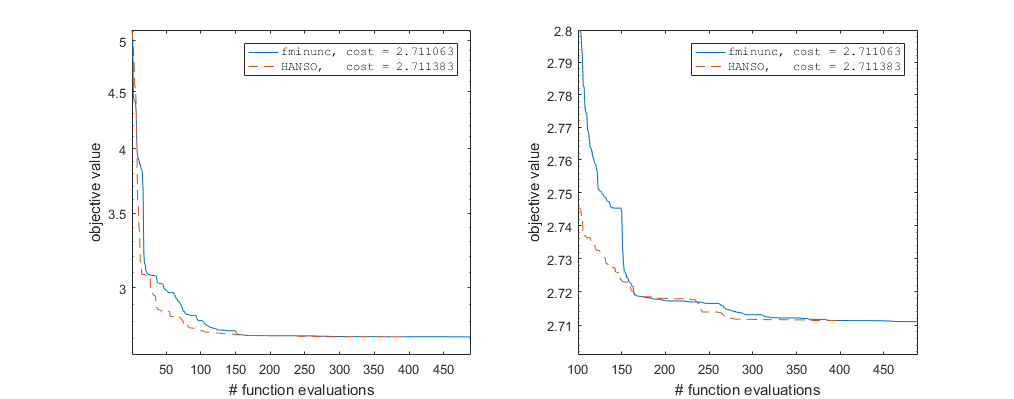}
\end{center}
\caption{Best objective value versus the number of function
  evaluations}\label{fig:h0025cost} 
\end{figure}

\begin{figure}[h]
\begin{center}
 \includegraphics[scale=0.6]{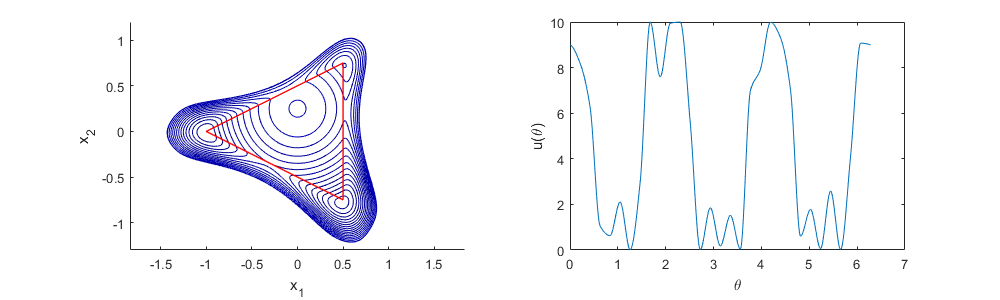}
\end{center}
\caption{Optimized domain (left) and control (right) computed using fminunc}\label{fig:h0025}
\end{figure}

\begin{figure}[h]
\begin{center}
 \includegraphics[scale=0.6]{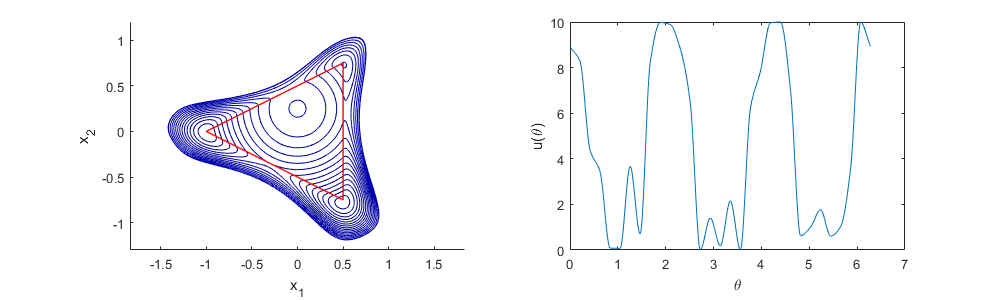}
\end{center}
\caption{Optimized domain (left) and control (right) computed using HANSO~2.2}\label{fig:h0025H}
\end{figure}

This time both optimizers end up to almost same final objective value
using essentially same amount of computational work. However, the
optimized controls and the corresponding domains $\Omega_\eps^\star$
differ considerably.

\end{example}

The results of Examples~1 and 2 show that the method works reasonably
well. The obtained shapes are  approximately feasible (due to
penalization) and the objective function value was substantially
reduced. However, there is no guarantee that any of the obtained
domains is even a local minimizer. 

From the results, we can also deduce that
the boundary control problem is ill-conditioned. The
effect of ill-conditioning can be seen by comparing the optimized
controls obtained by using the two different optimizers. Comparing the
convergence history we see that significant differences in the details
of the boundary control have only a minor effect on the value of the
objective function. This is due to the fact that the solution of the
Poisson problem is generally smoother than the data. Similar stiffness
appears e.g.\ in identification problem related to the Bernoulli free
boundary problem studied in \cite{THM08}.


\section{Conclusions}

In this paper we have considered some computational aspects of a shape
optimization problem with the 
state constraint given by a free boundary/obstacle problem. The free
boundary problem is formulated as a quadratic programming problem
which is then solved using the state-of-the-art tools. To solve the
optimal shape design problem, a boundary control approach was used.
Its main advantage is that there is no need to  consider moving
domains. This is advantageous especially from the computational point
of view. 

A well-known feature of optimal control problems governed by obstacle type
problems is that the control-to-state mapping is not smooth in
general. However, in discrete setting, it is piecewise smooth. In
those points where it is smooth, the gradient of the objective
function can be evaluated in a
straightforward way using the adjoint approach.

Numerical examples show that the location of the boundary of the
contact zone can be adjusted by changing the boundary
control. However, it seems that the problem is ill-conditioned in the
sense that relatively large changes in the boundary control have only
a little effect on the location of the contact zone boundary. 
We did not make a detailed trade-off study between accuracy and
oscillations. That will be a topic in further studies. 
One
possible remedy (and a topic for further studies) would be to consider
distributed control instead of boundary control as was also done in
\cite{NT95}. 

\clearpage

\bibliographystyle{siam}

\end{document}